\documentstyle[12pt, leqno,amscd,amsopn,upref]{amsart}
\def\qed{\relax
     \ifmmode
       ~\hfill\Box
     \else
        \unskip\nobreak ~\hfill$\Box$%
      \fi \par}

\theoremstyle{plain}
\newtheorem{prop}{Proposition}
\newtheorem{thm}{Theorem}
\newtheorem*{lemma}{Lemma}

\theoremstyle{definition}
\newtheorem{defn}{Definition}
\newtheorem{exmp}{Example}
\newtheorem*{rem}{Remark}

\newcommand{\neu}[1]{{\rm (\ref{#1})}}
\newcommand\half{\mbox{\small{$\frac{1}{2}$}}}
\newcommand\vts{\mkern1mu}
\newcommand\ts{\mkern 2mu}
\newcommand\End{{\sf {End}}\vts }
\newcommand\Hom{{\sf {Hom}}\vts }
\newcommand\Spin{{\sf {Spin}}\vts }

\newcommand\Cpin{{\mathbb C}^{\times}\cdot{\sf {Pin}}\vts }
\newcommand\Pin{{\sf {Pin}}\vts }
\newcommand\Pinc{{\sf {Pin}}^{c}\vts }
\newcommand\Lpin{{\sf {Lpin}}\vts }
\newcommand\SO{{\sf {SO}}\vts }
\newcommand\Ort{{\sf {O}}\vts }
\newcommand\Cl{{\sf {Cl}}\vts}
\newcommand\U{{\sf {U}}\vts}
\newcommand\GL{{\sf {GL}}\vts}
\newcommand\tbwedge{\mbox{\large $\wedge $}}
\newcommand{\Vol}{Vol}
\newcommand{\vol}{vol}
\newcommand{\Ad}{Ad}
\newcommand{\Obj}{Obj}
\newcommand{\Mor}{Mor}
\newcommand{\id}{id}
\newcommand{\Aut}{Aut}

\begin{document}

\noindent {\large \bf Spin spaces, Lipschitz groups, and spinor \\bundles}\\

\noindent Thomas Friedrich\\
{\small  \it Institut f\"ur Reine Mathematik, Humboldt Universit\"at\\
Ziegelstrasse 13A, 10099 Berlin,  Germany\\
(E-mail: friedric@@mathematik.hu-berlin.de)}\\

\noindent Andrzej Trautman\\
{\small \it Instytut Fizyki Teoretycznej, Uniwersytet Warszawski\\
Ho\.za 69, 00681 Warszawa, Poland\\
(E-mail: amt@@fuw.edu.pl)}\\

\noindent {\small {\bf Abstract.}  It is shown that every bundle
	 \(\varSigma\to M\) of complex spinor modules over the Clifford
	 bundle \( \Cl(g) \) of a Riemannian space \((M,g)\) with local
	 model \((V,h)\)
	is associated with an lpin (``Lipschitz") structure on $M$, this
	being a reduction of the \(\Ort(h)\)-bundle of all orthonormal
	frames on \(M\) to the Lipschitz group \(\Lpin(h)\) of all
	automorphisms of a suitably defined spin space.  An explicit
	construction is given of the total space of the
	\(\Lpin(h)\)-bundle defining such a structure.  If the dimension
	\( m \) of \( M \) is even, then the Lipschitz group coincides
	with the complex Clifford group and the lpin structure can be
	reduced to a  pin\(^{c}\) structure.  If \( m=2n-1 \), then a
	spinor module \( \varSigma \) on \( M \) is of the Cartan type:
	its fibres are \( 2^{n} \)-dimensional and decomposable at every
	point of \( M \), but the homomorphism of bundles of algebras \(
	\Cl(g)\to\End\varSigma \) globally decomposes if, and only if, \(
	M \) is orientable.  Examples of such bundles are given.  The
	topological condition for the existence of an lpin structure on an
	odd-dimensional Riemannian manifold is derived and illustrated by
	the example of a manifold admitting such a structure, but no pin\(
	^{c} \) structure.
\\

\noindent {\bf Keywords:} Clifford and spinor bundles, spin structures, pin\({}^{c }\)
structures\\

\noindent {\bf Mathematics Subject Classification (1991):} {Primary 15A66 and 53A50; Secondary 81R25 and
83C60.}\\}

%\end{opening}
\noindent {\em Dedicated to the memory of Alfred Gray}\\

%\vspace{0.5cm}

\section{Introduction}
\noindent There are at least two approaches to objects
 appearing in classical differential
geometry.  One of them, most common now,
can be traced back to the method
of ``moving frames'' of Darboux and Cartan; it is based on the notion
of a  {\it
principal\/} bundle with an infinitesimal connection in the sense of
Ehresmann. Fields of geometric objects are defined as sections of
associated bundles. Another approach assumes {\it vector\/} bundles
with
connections as the starting point; as
emphasized by Lang \cite{La}, this approach is
well adapted to the treatment of
infinite-dimensional manifolds, modeled on Banach spaces.  In the
category of finite-dimensional smooth manifolds, the two
approaches are essentially equivalent.

Spinors  on Riemannian manifolds can also be introduced in  two
similar
ways; however, the relation between the two approaches is subtler
than
in the case of tangent vectors. It is the purpose of this paper to
describe this relation in the general setting of
(not necessarily proper and orientable) Riemannian manifolds.
\smallskip

Recall first that if  an orientable, proper Riemannian \( m
\)-manifold
\( (M,g)\) has a spin structure (``if \( M\) is spin"),

\begin{equation}
\begin{CD}
\Spin_{m}      @>>>Q_{0}\\
 @V\Ad VV  @VV {\chi}_{0} V\\
 \SO_{m}       @>>>P_{0}    @>\pi >>M,
 \end{CD}
\label{e:spinstr}
\end{equation}
 then, given a spinor representation \( \theta \) of \( \Spin_{m}\) in
 a vector space \( S_{0}\), one defines spinor fields as sections of
 the associated vector bundle \( Q_{0}\times_{\theta} S_{0}\to M\).
 Another definition focuses on the vector bundle itself: it assumes
 the existence of a bundle \( \varSigma\to M\) of modules over the
 bundle \( \Cl(g)\) of Clifford algebras on \( M\).  This point of
 view can be traced back to early papers by physicists \cite{Fo}; in
 particular, to the remarkable article by Schr\"odinger \cite{Sch}
 which contains the first---to our knowledge---derivation of the
 formula for the square of the Dirac operator on Riemannian
manifolds.
 There are two somewhat forgotten papers by Karrer \cite{Kar1,Kar2} on
 the description of spinors in terms of Clifford modules; they contain
 the relevant definitions in the language of contemporary mathematics;
 see also Chapter 3 in \cite{BGV} and Chapter II in \cite{LM}.  The
 second approach is more general in the sense that the fibres of \(
 \varSigma\) need not be isomorphic to a spinor space, carrying an
 irreducible representation of the Clifford algebra; for example, one
 can take for \( \varSigma\) the bundle \( \tbwedge T^{*}M\) of
 exterior algebras.  This possibility was considered by physicists, in
 an attempt to find the relativistic, quantum-mechanical equation of
 the electron, as early as 1928 \cite{IL}.  Later, it has been
 developed by K\"ahler \cite{Kh}.  One easily sees that a bundle
 associated by a spinor representation with a spin structure is a
 bundle of modules over \( \Cl(g)\) (see, e.g., Prop.  3.8 in
 \cite{LM}), but the converse is not true, even if \( \Cl(g_{x})\to
 \End \varSigma_{x}\) is the spinor representation for every \( x\in
 M\) (see Examples  \ref{exmp:2},  \ref{ex:proj} and \ref{ex:grass}
 in  this  paper).

To compare these two definitions, it is convenient to introduce a {\it
category of spin spaces\/} (Section \ref{SpinS}).
We show that a
spinor bundle \( \varSigma\) of spinor modules over \( \Cl(g)\) is
associated with an lpin structure on \( M\), this being a reduction of
the \( \Ort_{m}\)-bundle of all orthonormal frames \( P\) to the {\it
Lipschitz group\/} \( \Lpin_{m} \) of automorphisms of a spin space.
If \( m \) is {\it even}, then \( \Lpin_{m} \) is essentially the {\it
Clifford group\/}\footnote{This name was introduced by Chevalley
\cite{Ch}.  R. Lipschitz was the first to consider groups associated
with Clifford algebras; see \cite{Weil} and the references given
there.  For this reason, we find it appropriate to associate his name
with one of the spinor groups.}.  The case of \( m \) {\it odd\/} is
somewhat more complicated because the complexification of the Clifford
algebra of a real, odd-dimensional vector space is not simple and the
adjoint representation of the corresponding Pin group does not cover
the full orthogonal group.  In this case, \( \Lpin_{m} \) is the
smallest group containing the Clifford group and such that the
reflection \( v\mapsto -v \) extends to an inner automorphism.
In Sections \ref{CBM} and \ref{SpBun} we recall the definition of a
bundle of Clifford modules on a Riemannian manifold and of a suitably
generalized spinor structure and spinor bundle. We show that every
faithful spinor bundle is associated with an lpin structure and prove
a theorem on the connection between the
orientability of \( M \) and the decomposability of such a spinor
bundle. An lpin structure on a manifold that is
 orientable or even-dimensional can be reduced to a
  pin\( ^{c} \) structure.  In Section \ref{Top} we establish
  the topological condition
for the existence of an lpin structure on an odd-dimensional,
non-orientable manifold, and give examples of such manifolds that
admit an lpin structure, but no pin\( ^{c} \) structure.

 \section{Notation and preliminaries}\label{s:Not}

\noindent We use a notation and terminology which are largely standard in
 differential geometry \cite{KN} and spinor analysis \cite{Fr2,LM}.
 If \( S\) and \( S'\) are finite-dimensional complex vector spaces,
 then \( \Hom (S,S')\) is the vector space of all complex-linear maps
 of \( S\) into \( S'\) and \( \End S=\Hom (S,S)\) is an algebra over
 \( {\mathbb C}\).  Every algebra under consideration here has \newpage \noindent a unit
 element; homomorphisms of algebras map units into units.
  We write \( S^{*}=\Hom(S, {\mathbb C}\ts)\); if \(
 f\in\Hom (S,S')\), then \( f^{*}\in \Hom ({S'}^{*},S^{*})\) is
 defined by \( \langle s,f^{*}(t')\rangle\) for every \( s\in S\) and
 \( t'\in {S'}^{*}\).  A similar notation is used for real vector
 spaces.

 A {\it quadratic space\/} is defined as a pair \( (V,h)\), where \(
 V\) is a real, finite-dimensional vector space and \( h:V\to \mathbb
 R\) is a non-degenerate quadratic form.  We denote by \( \hat{h}\)
 the symmetric linear isomorphism of \( V\) onto \( V^{*}\),
 associated with \( h\).  For our purposes, it is convenient to label
 with \( h \) the groups and algebras associated, in a natural manner,
 with the quadratic space \( (V,h)\); e.g.  \( \Ort(h) \) is the group
 of all its orthogonal automorphisms.  The real {\it Clifford
 algebra\/} \( \Cl(h)\) corresponding to \( (V, h)\) contains \(
 {\mathbb R}\oplus V\) as a vector subspace and \( v^{2}=h(v)\) for
 every \( v\in V\).  If \( (e_{1},\dots,e_{m}) \) is an orthonormal
 frame in \( V \), then \( \eta=e_{1}\dots e_{m}\in\Cl(h) \) and \(
 -\eta \) are the {\it volume elements}.  Their squares are either 1
 or \( -1 \), depending on the signature of \( h \); we define \(
 \iota(h)\in\{1,\sqrt{-1}\} \) so that \( \eta^{2}=\iota(h)^{2} \).
 If \( V ={\mathbb R}^{m} \) and \( h \) is the standard positive
 (resp., negative) definite quadratic form on that vector space, then
 we write \( \Cl^{+}_{m} \) (resp., \( \Cl^{-}_{m} \)) instead of \(
 \Cl(h)\) and use similar conventions for the various groups
 associated with \( (V,h) \); see \neu{e:spinstr} for examples.  The
 isometry \( v\mapsto -v\) extends to an involutive automorphism \(
 \alpha\) of the algebra defining its \( {\mathbb Z}_{2}\)-grading: \(
 \Cl(h)= \Cl^{0}(h)\oplus \Cl^{1}(h)\).  The even subalgebras of \(
 \Cl^{+}_{m} \) and \( \Cl^{-}_{m} \), which are isomorphic, are
 denoted by \( \Cl^{0}_{m} \).  The map \( {\mathbb R}^{m}\to
 \Cl^{0}_{m+1} \), \( v\mapsto ve_{m+1} \), extends to an isomorphism
 of algebras
 \begin{equation}
 i_{m}:\Cl^{-}_{m}\to\Cl^{0}_{m+1}.
 \label{e:ia}
 \end{equation}
 Let \( {\mathbb C}^{\times}={\mathbb C}\setminus \{0\} \) be the multiplicative group of complex
 numbers.  If \( G \) and \( H \) are groups and the sequences \(
 1\to{\mathbb Z}_{2}\to G \) and \(1\to{\mathbb Z}_{2}\to H \) are exact,
 then there is also the exact sequence \(1\to{\mathbb Z}_{2}\to G\cdot H
 \), where
 \[
 G\cdot H \quad\mbox{is the group}\quad (G\times
 H)/{\mathbb Z}_{2}.
\]
One  writes \( G^{c} \) instead of \( \U_{1}\cdot G \).

 All manifolds and maps among them are assumed to be smooth.
 Manifolds are finite-dimensional, but not necessarily compact.  If \(
 f:M\to N \) is a map of manifolds, then \( Tf:TM\to TN \) is the
 corresponding tangent (derived) map.  A {\it Riemannian manifold\/}
 \( M \) is assumed to be connected; it has a metric tensor field \( g
 \) which is non-degenerate, but not necessarily definite; if it is,
 then \( M \) is said to be {\it proper\/} Riemannian.  If \( \pi:E\to
 M\) is a fibre bundle over a manifold \( M\), then \( E_{x}=
 \pi^{-1}(x)\subset E\) is the fibre over \( x\in M\); in particular,
 \( T_{x}M\subset TM\) is the tangent vector space to \( M\) at \(
 x\).  If the pair \( (f,h) \) is a a morphism of the principal
 bundles \( G\to Q\to M \) and \( G'\to Q'\to M \) so that \( f:G\to
 G' \) is a morphism of Lie groups and \( h:Q\to Q' \) is a morphism
 of fibre bundles such that \( h(qa)=h(q)f(a) \) for every \( q\in Q
 \) and \( a\in G \), then one says that \( Q \) is an \(f\)-{\it
 reduction\/} of \( Q' \) to the Lie group \( G \).  Let \( G\to Q\to
 M \) be a principal bundle and let \( \theta \) be a representation
 of \( G \) in a complex vector space \( S \); the total space of the
 vector bundle associated by \( \theta \) with \( Q \) is the set \(
 Q\times_{\theta} S \) of all classes \( [(q,s)] \), where \( (q,s)\in
 Q\times S \) and the equivalence relation is: \((q,s)\equiv (q',s')\)
 if, and only if, there is \( a\in G \) such that \( q'=qa \) and \(
 s=\theta(a)s' \).

 To reconstruct the principal bundle from a vector bundle, it is
 convenient to introduce a category \( {\mathcal  C} \) of isomorphisms
 of finite-dimensional vector spaces such that if \( S\in \Obj\mathcal
 C \), then \( G= \Mor (S,S)\) is a closed subgroup of \( \GL(S) \).
 Given a vector bundle \(\varSigma\to M \) with typical fibre \(S \in
 \Obj {\mathcal  C} \), such that \( \varSigma_{x}\in \Obj {\mathcal  C} \)
 for every \( x\in M \), one defines the total space of the principal
 bundle to be
 \begin{equation}
	Q=\{q\in\Mor (S,\varSigma_{x})\vert x\in M\}.
	\label{e:defQ}
\end{equation}

The group \( G \) acts in \( Q \) by composition and there is a
natural projection \( Q\to M \).

Throughout this paper, given a positive integer \( m\), we define \(
\nu(m)\) to be the integer part of \( \half(m+1) \); sometimes we
write \( n \) instead of \( \nu(m) \).

\section{Spinor representations and spin spaces}\label{SpinS}

\noindent For the purposes of this paper, it is convenient to have a precise
definition of spinor groups and representations that, in the case of
odd dimensions, slightly difers from the one prevailing in the
literature.
\smallskip

 \noindent (i) Recall that if the dimension \( m\) of a quadratic
 space \(( V,h)\) is {\it even}, \( m=2n\), then the algebra \( \Cl
 (h) \) is central simple and has one, up to complex equivalence,
 faithful and irreducible {\it Dirac\/} representation in a complex,
 \( 2^{n}\)-dimensional vector space \( S\).  Restricted to \(
 \Cl^{0}(h) \), this representation decomposes into the direct sum of
 two complex-inequivalent, {\it half-spinor\/} (Weyl) representations
 \( \sigma_{+} \) and \( \sigma_{-} \).
 \smallskip

 \noindent (ii) If \( m \) is {\it odd}, \( m=2n-1 \), then the
 algebra \( \Cl^{0}(h) \) is central simple; it has a faithul and
 irreducible {\it Pauli\/} representation in a complex vector space \(
 S_{0} \) of dimension \( 2^{n-1} \); this representation extends to
 two complex-inequivalent representations \( \sigma \) and \(
 \sigma\circ\alpha \) of the full algebra \( \Cl(h) \) such that \(
 \sigma(\eta)=\iota(h)\id_{S_{0}} \).  These representations need not
 be faithful (example: \( \Cl^{+}_{1} \)).  The {\it Cartan\/}
 representation\footnote{This representation rarely appears because it
 is decomposable; it is needed to define the Dirac operator on
 odd-dimensional, non-orientable pin manifolds \cite{Ca,T}.  The names
 of Pauli and Dirac are associated by physicists with spinors in
 dimensions 3 and 4, respectively.} \( \sigma\oplus\sigma\circ\alpha
 \) of \( \Cl(h) \) in the \( 2^{n}\)-dimensional vector space \(
 S=S_{0}\oplus S_{0}\) is faithful.

The following Proposition is an obvious  consequence of
\neu{e:ia} and our terminology:
\begin{prop}\label{p:CDP}
Consider the sequence of homomorphisms of algebras
\[
\begin{CD}
\Cl^{-}_{m}@>{\mathrm {inj}}>>\Cl^{-}_{m+1}
@>{i_{m+1}}>>\Cl^{0}_{m+2}\stackrel{\theta}
{\longrightarrow}\End S.
\end{CD}
\]
If \( m \) is even (resp., odd) and \( \theta \) is a Weyl
(resp., Pauli)
representation, then \( \theta\circ i_{m+1} \) is a Pauli (resp.,
Dirac) representation and \( \theta\circ i_{m+1}\circ{\mathrm {inj}}
\) is a Dirac (resp., Cartan) representation.
\end{prop}

\begin{defn}\label{d:reprCl}
Any representation of \( \Cl(h) \) or \( \Cl^{0}(h) \) equivalent to
one of the representations described in (i) and (ii) is called a {\it
spinor representation\/} of that algebra.
\end{defn}

\begin{defn}
 Let \( (k,l) \) be a pair of non-negative integers and let \( m=k+l
 \).  The {\it category\/} \( {\mathcal  C}_{k,l} \) {\it of spin spaces}
 is as follows.  An {\it object\/} of \( {\mathcal  C}_{k,l}\) is a
 triple \( \varsigma =(S,V,h)\) such that \( S \) is a complex, \(
 2^{\nu(m)} \)-dimensional vector space, \( V\subset\End S\) is a real
 vector space of dimension \( m\) and \( v^{2}=h(v)\id_{S} \) for
 every \( v\in V \), where \( h \) is a quadratic form on \( V \) of
 signature \( (k,l) \).  {\it Morphisms\/} between two spin spaces \(
 \varsigma =(S,V,h)\) and \( \varsigma' =(S',V',h')\) of the same
 category are defined by
 \[
 \Mor(\varsigma,\varsigma')=
 \{a\in\Hom(S,S')\vert \mbox{\(a\) is invertible and
 \(aVa^{-1}=V'\)}\}.
\]
\end{defn}
If \( a\in\Mor (\varsigma,\varsigma')\), then the map \( V\to V' \),
given by \( v\mapsto ava^{-1} \), is an isometry: there is a forgetful
functor from the category of spin spaces to the corresponding category
of quadratic spaces.
If \( (S,V,h)\in\Obj{\mathcal  C}_{k,l} \),
then \( (S,\sqrt{-1}\ts V,-h )\in\Obj{\mathcal  C}_{l,k}\).
Given a quadratic space \( (V,h)  \) of even (resp., odd) dimension
\( m \), one constructs a spin space by considering a Dirac (resp.,
Cartan) representation of \( \Cl(h) \) in \( S \) and {\it
identifying\/}  \( V \)  with its image in \( \End S \).

\begin{prop} Let \( (S,V,h) \) be a spin space.
The dimension of the complex vector space
\[
A(h)=\{w\in\End S\vert wv+vw=0\mbox{ for every }v\in V\}	
 \]
 is \( 1 \)  or \( 2  \) depending on whether \( m \) is {\em even}
or {\em odd}.
\end{prop}

{\it Proof.} 
If \( m \) is even, then \( A(h) \) is spanned by  \(
\varGamma=\sqrt{-1}\ts\iota(h)\eta \). If \( m \) is odd, then it is
convenient to represent the elements of \( \End S \) in a block form,
corresponding to the decomposition \( S=S_{0}\oplus S_{0} \), so that
 \begin{equation}
 v=\left(\begin{array}{cc}\sigma(v)&0\\0&-\sigma(v)\end{array}\right)
 \in V,\quad
 \eta =\iota(h)\left(\begin{array}{cc}I&0\\0&-I\end{array}\right),
 	\label{e:v&eta}
 \end{equation}
 where \( \sigma \) is the representation defined in  (ii) above and
\( I=\id_{S_{0}} \). The space \( A(h) \) is then spanned by the pair
\( (\varGamma, \varGamma\eta )\), where
\[ \varGamma= \left( \begin{array}{cc}0&-I\\I&0\end{array}\right). \]\qed

  Since both the Dirac and the Cartan representations are faithful, it
  follows from the universal property of Clifford algebras that, given
  a spin space \( (S, V, h)\), one can {\it identify\/} \( \Cl (h) \)
  with the subalgebra of \( \End S\) generated, over the reals, by \(
  V\subset\End S\).  After this identification, one has
\begin{equation}
\alpha(a)=\varGamma^{-1}a\varGamma\quad\mbox{for every}\quad
a\in\Cl(h)\subset\End S,
\label{e:alpha}
\end{equation}
where, and in the sequel, \( \varGamma \) is an element of \( A(h) \)
such that \( \varGamma^{2}=-\id_{S} \).  If \( m \) is even, then,
over the complex numbers, the real vector space \( V\) generates \(
{\mathbb C}\otimes\Cl(h)=\End S\); if \( m \) is odd, then, over \(
{\mathbb C} \), the vector space \( V\oplus {\mathbb R}\varGamma\subset
\End S \) generates the algebra \( \End S \).  The involutive
automorphism \( a\mapsto \eta a\eta^{-1} \) of \( \End S \) defines a
\( {\mathbb Z}_{2} \)-grading of the algebra, \( \End S=\End^{0}
S\oplus\End^{1} S \), so that
\begin{equation}
 \End^{\varepsilon}S=
 \{a\in\End S\vert a\eta=(-1)^{\varepsilon}\eta a\},\;
\varepsilon=0,1.
 \label{e:grading}
\end{equation}

Let \( (S,V,h) \) be a spin space.  The group \( \Pin(h) \) can be now
defined as the subgroup of \( \GL(S) \) generated by all unit vectors,
i.e.  by all elements \( v\in V\subset \End S\) such that either \(
v^{2}=\id_{S} \) or \( v^{2}=-\id_{S} \).  The {\it twisted adjoint\/}
representation \( \widetilde{\Ad} \) of \( \Pin(h) \) in \( V \) is
defined by
\[
\widetilde{\Ad}(a)v=\alpha(a)va^{-1},\; a\in\Pin(h)\mbox{
and } v\in V.
\label{e:twist}
\]
The spin group is \( \Spin(h)=\Pin(h)\cap \Cl^{0}(h) \)
and \( \Spin_{0}(h) \) is the connected component of \( \Pin(h) \);
if the form \( h \) is definite, then \( \Spin(h)=\Spin_{0}(h) \).

The linear map \( V\to\End S \), \( v\mapsto \varGamma v \), has the
Clifford property, \( (\varGamma v)^{2}=h(v)\id_{S} \).  It extends to
a representation of the Clifford algebra
\begin{equation}
\gamma:\Cl(h)\to \End S\quad\mbox{such that \( \gamma\circ\alpha=
\Ad(\varGamma^{-1} )\circ\gamma \) }.
	\label{e:dfgamma}
\end{equation}
Since \( \varGamma v=(\id_{S}+\varGamma)v(\id_{S}+\varGamma)^{-1} \),
the representation \neu{e:dfgamma} is equivalent to the inclusion
representation \( \Cl(h)\to\End S \).
Note also that if the dimension of \( V \) is odd, then
\[
\gamma(v)=\left(\begin{array}{cc}0&\sigma(v)\\\sigma(v)&0\end{array}	
\right)	\label{e:frmgamma}
\]
whereas, in the inclusion representation, a vector is
represented by a ``block-diagonal'' matrix \neu{e:v&eta}.

Let \( \Pin_{\gamma}(h) \) be the image of \( \Pin(h) \) by the
monomorphism \neu{e:dfgamma} of algebras; by restriction, it gives
rise to the isomorphism of  groups
\begin{equation}
  \gamma:\Pin(h)\to\Pin_{\gamma}(h).
\label{e:iso}
\end{equation}
 The isomorphic groups \( \Pin(h) \) and \( \Pin_{\gamma}(h) \) are
 differently situated in \( \End S \) relative to \( V \).
\begin{prop}
 The groups \( \Pin(h) \) and \( \Pin_{\gamma}(h) \) provide
 equivalent extensions of \( \Ort(h) \) by \( {{\mathbb Z}}_{2} \): the
 diagram of group homomorphisms
\[
\begin{CD}
1@>>>{\mathbb Z}_{2}@>>>{\Pin(h)}@>{\widetilde{\Ad}} >>{\Ort(h)}@>>>1\\
&&@|@VV{\gamma} V @VVV\\
1@>>>{\mathbb Z}_{2}@>>>{\Pin_{\gamma}(h)}
@>{\Ad} >>{\Ort(h)}@>>>1
\end{CD}
\]
is commutative and its two horizontal sequences are exact.
\label{p:diag}
\end{prop}

{\it Proof.} 
The exactness of the upper sequence is classical \cite{ABS}.  A unit
vector \( u\) is in \( \Pin(h) \), \( \varGamma u\) is in
\(\Pin_{\gamma}(h) \) and the equalities
\(\widetilde{\Ad}(u)v=-uvu^{-1}=\varGamma uv(\varGamma u)^{-1}=
\bigl(\Ad\circ\gamma\bigr)(u)v \) complete the proof. \qed 

\medskip

Recall that, for \( m \) odd, the adjoint representation \( \Ad \)
maps \( \Pin(h) \) onto \( \SO(h) \) with a four-element kernel.  In
every dimension,
\[
 \Pin_{\gamma}(h)\cap\Cl^{0}(h)=\Spin(h).
 \]
The  {\it extensions\/} of \( \Ort(h) \) by \( {\mathbb Z}_{2} \)
corresponding to \( \Pin(h) \) and \( \Pin(-h) \) are {\it
inequivalent}, even in the case of \( h \) of neutral signature \(
(k,k) \), when the {\it groups\/} \( \Pin(h) \) and \( \Pin(-h) \) are
{\it isomorphic}.  Restricted to \( \Spin(h) \), the representations
\( \Ad \) and \( \widetilde{\Ad} \) coincide; the groups \( \Spin(h)
\) and \( \Spin(-h) \) are isomorphic and provide equivalent
extensions of \( \SO(h) \) by \( {\mathbb Z}_{2}\).

\begin{defn}
Let \( (S,V,h) \) be a spin space. A closed subgroup
\( G \) of \( \GL(S) \) is called a {\em spinor
group\/} if it contains the group \( \Spin_{0}(h) \) and is such that
\( aVa^{-1}=V\) for every \( a\in G \).
A representation of \( G \) equivalent to the evaluation
representation of \( G \) in \( S \)---or to its
subrepresentation--- is called a {\em spinor
representation\/} of this group.\label{d:sgroup}
\end{defn}
If the group \( G\) is contained in \(\Cl(h)\subset\End S \), then its
spinor representation is equivalent to the restriction to \( G \) of
the representation described in either (i) or (ii), depending on the
parity of \( m \).  The groups \( \Spin(h) \), \( \Pin(h) \), \(
\Pin(-h) \), \( \Pin_{\gamma}(h)\), \( \Pinc(h)=\U_{1}\cdot \Pin(h) \)
and \( {\mathbb C}^{\times}\cdot\Pin(h) \) are spinor groups.

\begin{defn}
The {\it Lipschitz group}  is  the ``largest'' spinor group,
\[
	\Lpin(h)=\{a\in\GL(S)\vert aVa^{-1}=V\}.
	\label{e:deflpin}
\]
\end{defn}
\begin{exmp}
Consider \( (S,V,h)\in {{\mathcal  C}}_{1,0} \) so that \( S= {\mathbb C}^{2} \) and \( V={\mathbb R} \). The injection \( V\to\End S \) is
given by \(t\mapsto\mbox{\scriptsize\( \left(\begin{array}{cc}
t&\phantom{-}0\\0&-t\end{array}\right)\)} \) and
\[
\Lpin_{1,0}=\{a\in\GL_{2}({\mathbb C})\ts\vert\ts \mbox{either \( a=
\mbox{\scriptsize\(\left(\begin{array}{cc}
\lambda &0\\0&\mu\end{array}\right)\)} \)
or
\(a= \mbox{\scriptsize\(\left(\begin{array}{cc} 0&
\lambda\\\mu &0\end{array}\right)\)} \)}\}.
\]
Let \( \psi:{\mathbb Z}_{2}\to \Aut ({\mathbb C}^{\times}\times
{\mathbb C}^{\times}) \) be the homomorphism given by
\[
 \psi(-1)(\lambda,\mu)=(\mu,\lambda).
  \]
   The group \( \Lpin_{1,0}
\) is isomorphic to the semi-direct product
\({\mathbb Z}_{2}\times_{\psi}({\mathbb C}^{\times}\times
{\mathbb C}^{\times})\).
\end{exmp}
\begin{prop}
The group
\[
K(h)=\{a\in\GL(S)\vert av=va\;\;\mbox{for every}\;\; v\in V\}
\]
is isomorphic to \( {\mathbb C}^{\times} \) or \( {\mathbb C}^{\times}
\times {\mathbb C}^{\times}   \) depending on whether \( m \) is even
or odd.
\end{prop}
{\it Proof.} 
If \( m \) is even, then \( K(h)={\mathbb C}^{\times}  \). If \( m
\) is odd, then \( K(h)={\mathbb C}^{\times}(\id_{S}+ \iota(h)\eta  )\times {\mathbb C}^{\times}(\id_{S}-\iota(h)\eta)  \).
\qed

\medskip

\begin{prop}
Let \( G \) be a spinor group. The homomorphism \( \Ad:G\to\Ort(h) \),
\( \Ad(a)v=ava^{-1} \),  where \( a\in G \) and \( v\in V \), is
surjective if, and only if,
\begin{equation}
\Pin_{\gamma}(h)\subset K(h)\cdot G
\label{e:Adsurj}
\end{equation}
\end{prop}
{\it Proof.} 
If \neu{e:Adsurj} holds, then, for every unit vector \( v\in V \),
there is \( k\in K(h) \) and \( a\in G \) such that \( \varGamma v=ka
\); therefore \( \Ad(a)=\Ad(\varGamma v) \) and \( \Ad(G) \) contains
the reflection in every hyperplane. Conversely, if  \( a\in G \) is
such that \( \Ad(a) \) is the reflection in the hyperplane orthogonal
to a unit vector \( v \), then \( \Ad(a^{-1})\circ \Ad(\varGamma
v)=\id_{V} \) so that \( a^{-1}\varGamma v\in K(h) \). \qed \medskip

%\\
In particular, the Lipschitz group satisfies \neu{e:Adsurj}; since
this group contains also all invertible vectors, the kernel of \(
\Ad:\Lpin(h)\to\Ort(h) \) is \( K(h) \).

If \( a\) is in a spinor group \( G \) and \( \eta \) is a volume
element, then \( a\eta a^{-1} \) is also a volume element.  Therefore,
\( a\eta a^{-1} \) equals either \( \eta \) or \( -\eta \) and the
grading \neu{e:grading} induces a \( {\mathbb Z}_{2} \)-grading of the
spinor group,
\begin{equation}
 G=G^{0}\cup G^{1},\mbox{ where }
G^{\varepsilon}=G\cap\End^{\varepsilon}S,\; \varepsilon=0,1.
\label{e:grgr}
\end{equation}
The real orthogonal group \( \Ort(h) \) is \( {\mathbb Z}_{2} \)-graded
by the determinant and the map \( \Ad:G\to \Ort(h) \) preserves the
grading.

If \( m \) is even and \( G\subset\Cl(h) \), then \(
G^{\varepsilon}=G\cap\Cl^{\varepsilon}(h) \).  For \( m \) odd, the
grading \neu{e:grgr} is {\it not\/} induced from that of the
Clifford algebra: an element of \( \Lpin(h) \) is odd if, and only if,
it is of the form
\begin{equation}
\left(\begin{array}{cc}0&\lambda\sigma(a)\\\mu\sigma(a)&0
\end{array}\right),
\mbox{ for some }\lambda,
\mu\in{\mathbb C}^{\times}\mbox{ and }
 a\in\Pin(h).
	\label{e:Aodd}
\end{equation}
In this grading, the elements \neu{e:v&eta} of \( \Lpin(h) \) are
{\it even}.

There is a complex quadratic form \( h_{0} \) on \( A(h) \) such that,
for every, \( w\in A(h) \), one has \( w^{2}=h_{0}(w)\id_{S} \).  If
\( a\in \Lpin(h) \) and \( w\in A(h) \), then \(\kappa(a)w= awa^{-1}
\) is also in \( A(h) \); this defines a homomorphism \( \kappa\) of
the Lipschitz group to the complex orthogonal group
\(\Ort(h_{0},{\mathbb C}) \).

\begin{thm}\label{t:strLpin}
Let \( (V,h) \) be a quadratic space of dimension \( m\).

\noindent { \rm (i)} If  \( m \) is {\em even}, then there is a
split exact sequence
 \[
 1\to{\mathbb C}^{\times}\cdot\Spin(h)\to\Lpin(h)\stackrel{\kappa}
\longrightarrow{\mathbb Z}_{2}\to 1	
 \]
and the group \( \Lpin(h) \) is isomorphic to \( {\mathbb
C}^{\times}\cdot\Pin(h) \).

\noindent {\rm (ii)} If \( m \) is {\em odd}, then there is the
exact sequence
 \begin{equation}
1\to{\mathbb C}^{\times}\cdot \Pin(h)\to\Lpin(h)\stackrel{\kappa}
\longrightarrow{\mathbb Z}_{2}
\times_{\varphi}{\mathbb C}^{\times}\to 1,	
  	\label{e:lpo}
  \end{equation}
where the homomorphism \( \varphi: Z_{2}\to \Aut {\mathbb C}^{\times}\),
defining the semi-direct product group structure, is given by \(
\varphi(-1)z=z^{-1} \), \( z\in {\mathbb C}^{\times} \).  There is an
isomorphism of groups
\begin{equation}
\bigl(\Pin(h)\times_{\psi}({\mathbb C}^{\times}\times {\mathbb
C}^{\times})\bigr)/{\mathbb Z}_{2}\to\Lpin(h)
\label{e:pin2lpin}
\end{equation}
such that
\[
\mbox{if \( a\in\Pin(h)\setminus\Spin(h) \),
 then}\;\; \psi(a)(\lambda,\mu)=(\mu,\lambda).
\]
  \end{thm}

{\it Proof.} If \( a\in\ker \kappa\subset\Lpin(h) \), then \( a \) is even with
respect to the grading \neu{e:grgr} of \(G= \Lpin(h) \)
and commutes with \( \varGamma \).

\noindent (i) For \(
m \) even, the group \( \Ort(h_{0},{\mathbb C}) \) is \(
\Ort_{1}({\mathbb C})={\mathbb Z}_{2} \), the map \( \kappa \) is the
grading homomorphism and \( \Lpin(h)=\Cpin(h) \). Let \( u\in V \)
be a unit  vector, put \( v=u \) if \( h(u)= 1 \) and \(
v=\sqrt{-1}\ts u \) if
\( h(u)=-1 \). The map \( f:{\mathbb Z}_{2}\to \Lpin(h) \) such that
\(f(1)=1 \) and \( f(-1)=v \) is a splitting homomorphism.

\noindent (ii) Let now \( m \) be odd.
The set of all odd elements of \( \Lpin(h) \) generates the group;
every odd element of \( \Lpin(h) \)  is of the
form \neu{e:Aodd}.
The group \( \Ort(h_{0},{\mathbb C}) \) is now isomorphic to \(
\Ort_{2}({\mathbb C})\) and can be identified with the semi-direct
product \({\mathbb Z}_{2}\times_{\varphi}{\mathbb C}^{\times} \).
Under this identification, the homomorphism
 \(\kappa  \) maps \neu{e:Aodd} to \( (-1,
\lambda\mu^{-1}) \).
The injection of \( \Cpin(h)\) into \(\Lpin(h) \) is given by
\[
[(\lambda, a)]\mapsto\left(\begin{array}{cc}\lambda
\sigma(a)&0\\0&\lambda \sigma(a)
\end{array}\right)\mbox{ for \( \lambda\in{\mathbb C}^{\times} \) and \(
a\in\Pin(h) \)}.
\]
For \( a\in\Pin(h) \) odd, the isomorphism
\neu{e:pin2lpin} sends
\[
 [(a,\lambda,\mu)]= [(-a,-\lambda,-\mu)]
 \]
to \neu{e:Aodd}.	\qed \medskip

\section{Bundles of Clifford modules and spinor structures}
\label{CBM}

\subsection{Clifford bundles and modules}

\mbox{}\\

\noindent Let \( (M,g) \) be a  Riemannian manifold and let \( g_{x} \) be the
restriction of \( g \) to the vector space \( T_{x}M \). The Clifford
algebra associated with the quadratic space \( (T_{x}M,g_{x}) \) is
\( \Cl(g_{x}) \) and \( \Cl(g)=\bigcup_{x}\Cl(g_{x}) \) is the total
space of the {\it Clifford bundle\/} of \( (M,g) \) \cite{ABS}.
\begin{defn}
A  bundle of {\em Clifford modules\/} on the
 Riemannian space \( (M,g) \) is a
complex vector bundle \( \varSigma\to M \) with
a homomorphism of
bundles of algebras
\begin{equation}
	\tau:\Cl(g)\to\End\varSigma.
	\label{e:tau}
\end{equation}
\end{defn}
In other words, for every \( x\in M \), the vector space \(
\varSigma_{x} \) is a left module over the algebra \( \Cl(g_{x}) \).
Restricted to \( TM\subset\Cl(g) \), the map \( \tau \) is a {\it
Clifford morphism\/}, i.e.  a homomorphism of vector bundles such that
\[
\tau(v)^{2}=g_{x}(v)\id_{\varSigma_{x}} \quad\mbox{for every \(x\in
M \) and \(v\in T_{x}M\)}.
\]
It follows from the universal property of Clifford algebras that,
conversely, given a vector bundle \( \varSigma \) over \( M \) and a
Clifford morphism \( TM\to\End\varSigma \), one can extend it to a
homomorphism \neu{e:tau} of bundles of algebras.
\smallskip

The following examples are well known:
 \begin{exmp} \label{exmp:1}
 The bundle of {\it exterior algebras\/} on \( M\). Put
\( \varSigma =\tbwedge T^{*}M\) and define \( \tau\) by
\( \tau(v)\omega=v\lrcorner
\ts\omega+\hat{g}_{x}(v)\wedge\omega\) for \( v\in T_{x}M\)
and \( \omega\in\varSigma_{x}\).
\end{exmp}

\begin{exmp}\label{exmp:2}
 Let \( (M,g)\) be an {\it almost Hermitean
space\/} and let \( J\) be the associated orthogonal almost complex
structure.  Define \( N=\{n\in {\mathbb C}\otimes TM\vert J(n)=
\sqrt{-1}\ts n\}\) and put \( \varSigma=\tbwedge N\) .  The map \(
\tau\) given, for every \( n\in N\) and \( \omega\in \varSigma\), by
\( \tau(n+\bar{n})\omega=\sqrt{2}\ts
(\hat{g}(\bar{n})\lrcorner\ts\omega+n\wedge\omega)\) is a Clifford
morphism.
\end{exmp}

\begin{exmp}\label{exmp:3}
Let \((M,g)\) be a Riemannian manifold of dimension \( m \) with \(
(V,h) \) as its local model.  Consider a pin structure on \( M \)
defined, as in \cite{Ba}, to be an \( \widetilde{\Ad} \)-reduction of
the bundle \( P \) of all orthonormal frames on \( M \) to the group
\( \Pin(h) \),
\begin{equation}
\begin{CD}
	\Pin(h)@>>> \widetilde{Q}\\
	@V\widetilde{\Ad} VV @VV\widetilde{\chi} V\\
	\Ort(h)@>>>P@>\pi >>M,
\end{CD}
\label{e:defpinstr}
\end{equation}
so that \( \pi\circ \widetilde{\chi}:\widetilde{Q}\to M \) is a principal \(
\Pin(h) \)-bundle and \(
\widetilde{\chi}(qa)=\widetilde{\chi}(q)\circ\widetilde{\Ad}(a) \) for every
\( q\in\widetilde{Q} \) and \( a\in\Pin(h) \).  Let \(
\theta:\Cl(h)\to\End S \) be either the Dirac (\( m \) even) or the
Cartan (\( m \) odd) evaluation representation.  The associated bundle
\( \varSigma=\widetilde{Q}\times_{\theta}S \) is a bundle of Clifford
modules.  The Clifford morphism is defined as follows: if \( v\in
T_{x}M, q\in \widetilde{Q}_{x} \) and \( s\in S \), then
\begin{equation}
\tau(v)[(q,s)]=[(q,\varGamma \widetilde{\chi}(q)^{-1}(v)s)].
\label{e:deftau}
\end{equation}
To check that \( \tau \) is correctly defined by \neu{e:deftau},
take \( a\in\Pin(h)\subset\GL(S) \) and use \neu{e:alpha} to compute
\begin{eqnarray*}
\tau(v)[(qa,a^{-1}s)]&=&[(qa,\varGamma\widetilde{\chi}(qa)^{-1}(v)a^{-1}s)]
\\
&=&[(qa,\varGamma\widetilde{\Ad}(a^{-1})
(\widetilde{\chi}(q)^{-1}(v))a^{-1}s)]\\
&=&\tau(v)[(q,s)].
\end{eqnarray*}
\end{exmp}
 Note that \( \varGamma \) in \neu{e:deftau} is essential to undo
the twisting implied by the use of \( \widetilde{\Ad} \).

\subsection{Spinor structures}\label{s:sstr}

\mbox{}\\

\noindent Consider a Riemannian, not necessarily orientable, manifold \( (M,g)
\) and a spin space \( (S,V,h) \) such that \( (V,h)\) is a local
model of the Riemannian manifold.  Let \( G \) be a spinor group in
the sense of Definition \ref{d:sgroup} such that \(
\Pin_{\gamma}(h)\subset K(h)\cdot G \) and let \( P \) be the bundle
of all orthonormal frames on \( M \).
\begin{defn}
A {\em spinor \( G \)-structure\/} on \( (M,g) \)
is an \( \Ad \)-reduction \( Q \) of
the \( {\Ort}(h) \)-bundle \( P \) to the group \( G \),
\begin{equation}
\begin{CD}
	G@>>> Q\\
	@V\Ad VV @VV\chi V\\
	\Ort(h)@>>>P@>\pi >>M,
\end{CD}	
	\label{e:Gstr}
\end{equation}
so that  \( \pi\circ\chi :Q\to M\) is
 a principal \(G\)-bundle and \( \chi(qa)=\chi(q)\circ\Ad(a) \)
 for every \( q\in Q \) and \( a\in G \).\label{d:Gstr}
\end{defn}
For \( G= \Pin_{\gamma}(h) \), \( \Pin_{\gamma}^c(h) \) and \(
\Lpin(h) \) one shortens the expression ``spinor \( G \)-structure''
to pin, pin\(^c \) and lpin structure, respectively.
A pin structure on an
non-orientable Riemannian manifold is often defined
as in Example \ref{exmp:3}, as an  \( \widetilde{\Ad} \)-reduction
\( \widetilde{Q} \) of \( P \) to the group
\( \Pin(h) \), where \( \widetilde{\Ad} \) is
the {\it twisted\/} adjoint representation. That definition is
equivalent to ours, given by an  \( \Ad \)-reduction of
\( P \) to \( \Pin_{\gamma}(h) \). To show this explicitly, let us
 consider the  pin structure \neu{e:defpinstr}. We construct
a pin structure
in the sense of Definition \ref{d:Gstr}
for \( G=\Pin_{\gamma}(h) \), described by
the diagram
\[
\begin{CD}
	\Pin_{\gamma}(h)@>>> Q\\
	@V{\Ad} VV @VV\chi V\\
	\Ort(h)@>>>P@>\pi >>M,
\end{CD}
\]
by defining \( Q \) to be the bundle associated
with the bundle \( \widetilde{Q} \) by the isomorphism \neu{e:iso} so
that
 \[
 Q=\widetilde{Q}\times_{\gamma}\Pin_{\gamma}(h)\quad\mbox{and}\quad
 \chi([(q,\gamma(a))])=\widetilde{\chi}(qa)
 \]
  for  \( q\in Q\)  and  \(a\in\Pin(h)\).
 A proof of the equivalence of
these two definitions is based on Proposition \ref{p:diag}. In view
of these observations, we restrict ourselves to spinor structures
defined in terms of the (untwisted) adjoint representation, even in
the case of non-orientable, odd-dimensional Riemannian manifolds.

On an orientable Riemannian manifold, after reducing \( P \) to the
group \( \SO(h) \), one can consider spinor \( G \)-structures such
that \( \Ad(G)=\SO(h) \). For \( G=\Spin(h) \)
and \( \Spin^{c}(h) \), one
obtains the usual notion of spin and spin\( ^{c} \) structure,
respectively.

\section{Spinor bundles}\label{SpBun}

\begin{defn}
A bundle of Clifford modules  \neu{e:tau}  on
 the Riemannian manifold \( (M,g) \)
is said to be a {\em spinor bundle\/} over \( M \) if, for every \(
x\in M \), the restriction  of \( \tau  \) to
\( \Cl(g_{x}) \) is equivalent to a spinor representation
of the algebra in the sense of   Definition \ref{d:reprCl}. If,
moreover, \( \tau \) is injective, then the spinor bundle is said to
be {\em faithful}.
\end{defn}
\begin{exmp}\label{ex:sph}
Let \( {\mathbb R}^{m+1} \) be embedded in \( \Cl^{+}_{m+1} \) so that
\( {\mathbb S}_{m}=\{x\in {\mathbb R}^{m+1}\vert x^{2}=1\}\) and \(
T {\mathbb S}_{m}=\{(x,y)\in {\mathbb S}_{m}\times {\mathbb R}^{m+1}\vert
xy+yx=0\} \).  The Riemannian metric \( g_{m} \) on the \( m \)-sphere
is defined by \( g_{m}(x,y)=y^{2} \).  Depending on whether \( m \) is
even or odd, take \( \theta:\Cl^{0}_{m+1}\to {\End} S_{m} \) to be either
the Pauli or the Dirac representation so that \( \dim_{\mathbb C}S_{m}=2^{\nu(m)} \).  Therefore, for \( m \) even (resp., odd) and
every \( x\in {\mathbb S}_{m} \), the Clifford map \( {\mathbb R}^{m}\to {\End S}_{m} \) given by \( y\mapsto\sqrt{-1}\theta(xy) \),
where \( xy+yx=0 \), extends to a Dirac (resp., Cartan) representation
of \( \Cl^{+}_{m} \) in \( S_{m} \).  The trivial vector bundle \(
\varSigma_{m}= {\mathbb S}_{m}\times S_{m} \) is made into a faithful
spinor bundle by defining \( \tau_{m}: \Cl(g_{m})\to\End
\varSigma_{m}\) so that \( \tau_{m}(x,y)(x,s)=(x,\sqrt{-1}\theta(xy)s)
\).
\end{exmp}
The bundles of Clifford modules described in Examples \ref{exmp:2}
and \ref{exmp:3} are also spinor bundles. Since there are Hermitean
manifolds that are not spin (e.g. the even-dimensional complex
projective spaces), Example \ref{exmp:2} shows that there are spinor
bundles on orientable manifolds that are not
associated with a spin structure.
\begin{defn}\label{d:sph}
Two spinor bundles
\[
\tau_{1}:\Cl(g)\to\End \varSigma_{1} \;\;\mbox{and}\;\;
\tau_{2}:\Cl(g)\to\End \varSigma_{2}
\]
on \( M \) are said to be {\em
equivalent\/} if there is an isomorphism of vector bundles \( i:
\varSigma_{1} \to \varSigma_{2}\) intertwining \( \tau_{1} \) and \(
\tau_{2} \) so that, for every \( x\in M \), \( s\in \varSigma_{x} \)
and \( u\in T_{x}M \) one has \( i(\tau_{1}(u)s)=\tau_{2}(u)i(s) \).
\end{defn}

\begin{exmp}\label{ex:proj}
Let \( {\mathbb P}_{m}={\mathbb S}_{m}/{\mathbb Z}_{2} \) be the real
projective \( m \)-space.  The action of \( {\mathbb Z}_{2} \) on \(
{\mathbb S}_{m} \) lifts to its tangent bundle and \( T{\mathbb P}_{m} \)
can be identified with \( T{\mathbb S}_{m}/{\mathbb Z}_{2} \).  The metric
on the sphere descends to the corresponding projective space.  If \(
\theta \) is one of the representations referred to in Example \ref
{ex:sph} and now \( \varSigma_{m}={\mathbb P}_{m}\times S_{m} \), then
the formula \(
\tau^{\pm}_{m}([(x,y)])([x],s)=([x],\pm\sqrt{-1}\theta(xy)s) \)
defines on \( \varSigma_{m} \) two structures of faithful spinor
bundles over \( {\mathbb P}_{m} \).  The spinor bundles \( \tau^{+}_{m}
\) and \( \tau^{-}_{m} \) are inequivalent.
\end{exmp}
Since for \(m\equiv 1\bmod 4 \), \( m>1 \), the space \( {\mathbb P}_{m}
\) has no spin structure, the above construction provides another
example of a spinor bundle that is not associated with a spin
structure.  This example can be generalized to covering spaces with a
finite cyclic group of deck transformations: let \(M\) be an
odd-dimensional Riemannian manifold with the fundamental group
isomorphic to \({\mathbb Z}_p\).  If the universal covering space of \(
M \) admits a spin structure, then there exists a spinor bundle over
\( M \).
\begin{thm}
Let \( (M, g) \) be a Riemannian manifold with an lpin structure
\neu{e:Gstr}, \( G=\Lpin(h)\subset\GL(S)\), and let \(
\theta:\Lpin(h) \to\GL(S)\) be the evaluation representation.  There
is a natural structure of a spinor bundle on the associated vector
bundle \( \varSigma=Q\times_{\theta} S\to M \).
\end{thm}
{\it Proof.} 
The Clifford morphism \( \tau:TM\to\End \varSigma \) is defined by
\[
\tau(v)[(q,s)]=[(q,\chi(q)^{-1}(v)s)],
 \]
where \( v\in T_{x}M \), \(
q\in Q_{x} \) and \( s\in S \). \qed \medskip

The following theorem shows that, conversely, every spinor bundle can
be so obtained.
\begin{thm}
Every faithful spinor bundle \( \varSigma\to M\) on a
Riemannian manifold \( (M,g)\) with local model \( (V,h)\)
is isomorphic to the  bundle  associated, by the spinor
representation, with an lpin structure
on that manifold.
\end{thm}
{\it Proof.} Since the bundle is assumed to be faithful, \( \tau \) is injective
and one can identify \( TM \) with its image by \( \tau \) in \( \End
\varSigma \).  One then constructs the total space \( Q\) of an lpin
structure by taking for the fibre \( Q_{x}\) the set of all
isomorphisms of the spin space \( (S,V,h)\) onto the spin space \(
(\varSigma_{x},T_{x}M,g_{x})\).  The map \( \chi :Q_{x}\to P_{x}\) is
given by \( \chi(q):V\to T_{x}M\), \( \chi(q)v=qvq^{-1}\).  If \( q\)
and \( q'\in Q_{x}\), then \( q^{-1}q'\in \Lpin(h)\); the group \(
\Lpin(h)\) acts freely and transitively on \( Q_{x}\) and \(
\chi(qa)=\chi(q)\circ\Ad(a)\) for every \( a\in \Lpin(h) \).  It
remains to check that the associated bundle of spinors \(
Q\times_{\theta} S\to M\) is isomorphic to \( \varSigma\to M\): such
an isomorphism is given by \( [(q,s)]\mapsto q(s)\), where \( q\in Q\)
and \( s\in S\). \qed \medskip

For every spinor bundle \neu{e:tau} and \( x\in M\), the restriction
\(\tau_{x}=\tau\vert \Cl(g_{x})\) is a spinor representation of the
Clifford algebra \( \Cl(g_{x})\) in the vector space \( \varSigma_{x}
\).  Therefore, if the dimension \( m \) of \( M \) is {\it odd}, then
this representation decomposes into two irreducibles.  Similarly, if
\( m \) is {\it even}, then the restriction \( \tau_{x}^{0} \) of \(
\tau_{x} \) to the even subalgebra \( \Cl^{0}(g_{x}) \) decomposes
into the sum of two half-spinor representations.  Put \( \Cl^{0}(g)
=\bigcup_{x}\Cl^{0}(g_{x})\) and define
\begin{equation}
	\Cl^{\mathrm{d}}(g)=\left\{
	\begin{array}{ll} \Cl(g) &\mbox{if \( m \) is odd,}\\
	\Cl^{0}(g)&\mbox{if \( m \) is even}.\end{array}\right.
\label{e:v}
\end{equation}
\begin{defn}\label{d:decomp}
Let \( \varSigma \) and \( \mathcal A\) be a vector bundle and a
bundle of algebras over \( M \), respectively. Let \(
\tau:\mathcal A\to\End\varSigma  \)
be a morphism of bundles of algebras. We say that \( \tau \) is {\em
decomposable\/} if there are two vector bundles \( \varSigma_{+} \)
and \( \varSigma_{-}  \) of positive fibre dimensions such that
\( \varSigma=\varSigma_{+}\oplus \varSigma_{-} \) and, for every \(
x\in M \), one has \( \tau({\mathcal A}_{x})\varSigma_{\pm x}\subset
\varSigma_{\pm x} \).
\end{defn}
\begin{thm}
Let \neu{e:tau} be a spinor bundle on an \( m \)-dimensional
Riemannian manifold \( M \) and let \( \tau^{\mathrm{d}} \) be the
restriction of \( \tau \) to the bundle of algebras \(
\Cl^{\mathrm{d}}(g) \) as in  \neu{e:v}.
A necessary and sufficient condition for \( \tau^{\mathrm{d}} \) to
be {\em decomposable\/} is that \( M \) be {\em orientable\/}.
\end{thm}

{\it Proof.} Let the quadratic space \( (V,h) \) be a local model of \( (M,g) \).
 If \( M \) is orientable, then there is a volume map \( \vol:M\to
\Cl^{\mathrm{d}}(g) \) such that \( \vol(x)^{2}=
\iota(h)^{2} \) for every \(
x\in M \). Defining
\begin{equation}
 \varSigma_{\pm}=(\id\pm \iota(h)\vol )\varSigma\quad\mbox{so that}
 \quad  \varSigma=\varSigma_{+}\oplus\varSigma_{-}
 \label{e:decom}
\end{equation}
and noting that \( \vol(x) \) is in the center of \(
\Cl^{\mathrm{d}}(g_{x}) \),
one obtains the required decomposition.
Conversely, assume that \( M \) is not orientable. Let
\(\Vol\subset  \Cl^{\mathrm{d}}(g)  \) be the
set of  normalized  volume elements
at all points of \( M \); the map \( \Vol\to M \) is a
double cover.  Non-orientability of \( M \) is
equivalent to the statement that the set \(\Vol\)
is connected.
 If \( \tau^{\mathrm{d}}  \) decomposes,  \( \tau^{\mathrm{d}} =
 \tau_{+}^{\mathrm{d}}\oplus \tau_{-}^{\mathrm{d}}\) and
 \( \varSigma_{+ }=
 \tau_{+}^{\mathrm{d}}(\Cl^{\mathrm{d}}(g) )\varSigma\),
  then
 \( \tau_{+}^{\mathrm{d}}(\Vol_{x}) =\{\iota(h)\id_{\varSigma_{+ x}},
 -\iota(h)\id_{\varSigma_{+ x}} \}\) so that  \(
 \tau_{+}^{\mathrm{d}}\) maps the connected set \( \Vol \) onto a set
 with two components. The contradiction shows that
 \( \tau^{\mathrm{d}}  \) does not decompose. \qed \medskip

\begin{rem}
A Cartan spinor bundle  \( \varSigma \to M\)
 associated with a pin
structure on an odd-dimensional manifold \( M \) always admits a
decomposition into two Pauli bundles, \( \varSigma=
\varSigma'_{+}\oplus \varSigma'_{-}\), corresponding to the
decomposition of the Cartan representation, \(
\gamma=\sigma\oplus\sigma\circ\alpha \). This decomposition of \(
\varSigma \)
is given by
\( [(q,s)]=[(q,s'_{+})]+[(q,s'_{-})] \), where \(
s'_{\pm}=\half(\id_{S}\pm\iota(h)\gamma(\eta))s \). If \( v\in T_{x}M
\), then \( \tau(v) \) maps \( \varSigma'_{+x} \)
 into \( \varSigma'_{-x} \).
Therefore,  on an orientable odd-dimensional \( M \),
this decomposition is different
from (``transversal to'') the decomposition \neu{e:decom}. The
following example illustrates this remark.
\end{rem}
\begin{exmp}
Let \( M \) be the real projective quadric \( ({\mathbb S}_{1}\times{\mathbb S}_{2})/{\mathbb Z}_{2} \) with a proper Riemannian
metric descending from \( {\mathbb S}_{1}\times{\mathbb S}_{2} \). This
is a non-orientable 3-manifold with a pin structure \cite{Ca}. A
typical
element of \( TM \) can be written as \( [(x,y,\xi,\eta)] \), where
\[
x,\xi\in{\mathbb R}^{2}\subset\Cl^{+}_{2},\; y,\eta\in {\mathbb R}^{3}\subset\Cl^{+}_{3},\; x\xi+\xi x=0,\; y\eta+\eta y=0,
\]
and \([(x,y,\xi,\eta)]=[(-x,-y,-\xi,-\eta)]  \). Let
\(\theta:\Cl^{+}_{2}\to\End S_{1} \) and \( \sigma:\Cl^{0}_{3}\to\End
S_{2} \) be, respectively, the  Dirac
and the Pauli representations in the
complex, 2-dimensional spaces of spinors \( S_{1} \) and \( S_{2} \).
Let \( (e_{1},e_{2},e_{3}) \) and \( (f_{1},f_{2}) \) be \
orthonormal
bases in \( {\mathbb R}^{3} \) and \( {\mathbb R}^{2} \), respectively.
 There is a trivial spinor bundle \( \varSigma=M\times
(S_{1}\otimes S_{2}) \) on \( M \) such that
\begin{eqnarray*}
\lefteqn{\tau([(x,y,\xi,\eta)])([(x,y)], s_{1}\otimes s_{2})}\\
& &=([(x,y)],
\theta(\xi)s_{1}\otimes\sigma(ye_{1}e_{2}e_{3})s_{2}+
s_{1}\otimes\sigma(y\eta)s_{2}).
\end{eqnarray*}
The map \( \varpi:M\to\End \varSigma\) defined by
 \[
 \varpi([(x,y)])([(x,y)],s_{1}\otimes
s_{2})=([(x,y)],\theta(xf_{1}f_{2})s_{1}\otimes
\sigma(ye_{1}e_{2}e_{3})s_{2})
\]
is involutive and \( \varpi([(x,y)])  \) anticommutes
 with \( \tau([(x,y,\xi,\eta)]) \).
Putting  \( \varSigma_{\pm}=
(\id_{\varSigma}\pm\varpi )\varSigma\),
one  obtains the decomposition referred to in the Remark. On the other
hand, the spinor bundle \( \varSigma \) is not decomposable in the
sense of Definition \ref{d:decomp}.
\end{exmp}

\section{Topological conditions}\label{Top}
\noindent In this section we restrict our considerations to {\it proper\/}
Riemannian
spaces. Following the notation of Section \ref{s:Not}, we write
\( \Pin^{+}_{m}=\Pin_{m,0} \) and \( \Pin^{-}_{m}=\Pin_{0,m} \).
These two groups provide inequivalent extensions of \( \Ort_{m} \) by
\( {\mathbb Z}_{2} \). It is known
that the groups \( \U_{1}\cdot \Pin^{+}_{m} \) and
\( \U_{1}\cdot \Pin^{-}_{m} \) are isomorphic and give equivalent
extensions of \( \Ort_{m} \) by \( \U_{1} \); therefore, it is
legitimate to denote both of these groups by \( \Pin_{m}^{c} \).
Similarly, the isomorphic groups \( \Lpin_{m,0} \) and \( \Lpin_{0,m}
\)
are denoted by \( \Lpin_{m} \). It follows from the isomorphism
\neu{e:pin2lpin} that \( \Lpin_{m} \) contains \(
(\Pin_{m}\times_{\psi}(\U_{1}\times U_{1}))/{\mathbb Z}_{2} \) as its
maximal compact subgroup.

If \( E\to M \) is a real vector bundle, then \( w_{i}(E)\in
{\mathrm H}^{i}(M,{\mathbb Z}_{2}) \) denotes its \( i \)th
Stiefel-Whitney
class.
The manifold \( M \) is orientable if, and only if,  \( w_{1}(TM)=0
\).

\begin{lemma}
For every real  vector bundle \(E\) over a manifold \( M \)
the class
\(w_1(E)^2\in {\mathrm H}^2(M, {\mathbb Z}_2)\) is the {\rm mod 2}
reduction
of an element of \( {\mathrm H}^2(M,{\mathbb Z})\).
 \end{lemma}
{\it Proof.} 
Assume the bundle \( E \) to have  \( m \)-dimensional fibres
 and consider the line bundle \(F=\tbwedge^m E \).
 The direct sum  \( F\oplus F \) is an orientable
bundle of
fibre dimension 2. Its Euler class \( e(F\oplus F) \) is an integral
cohomology class and \( w_{2}(F\oplus F) \) is the mod 2 reduction of
\(e(F\oplus F) \). The
Whitney
product theorem gives \( w_{2}(F\oplus F)= w_{1}(F)^{2} \). The
equality \( w_{1}(E)=w_{1}(F) \) is established by the ``splitting
method'' used in \cite{Hi} in the proof of Theorem 4.4.3
for the Chern classes. \qed \medskip 

There are  well-known topological obstructions to
the existence of the various pin and spin
structures on a a Riemannian space. Recall that the necessary and
sufficient
conditions for the existence of these structure are as follows
\cite{Kb}:

\begin{enumerate}
\item[\ \ (i)] spin structure: \(w_1(TM)=0\) and  \(w_2(TM)=0\);\\
\item[\ (ii)] pin\( ^{+} \) structure: \(w_2(TM)=0\);\\
\item[(iii)] pin\( ^{-} \) structure: \( w_{1}(TM)^{2}+w_2(TM)=0\);\\
\item[\ (iv)] spin\({}^c\) structure: \(w_1(TM)=0\) and there exists a
cohomology class
\(c\in {\mathrm H}^2(M, {\mathbb Z})\) such that \(w_2(TM)\equiv c \bmod
2\);\\
\item[\ \ (v)]  pin\({}^c\) structure:  there exists a
cohomology class
\(c\in {\mathrm H}^2(M, {\mathbb Z})\) such that
\(w_2(TM)\equiv c \bmod 2\).
\end{enumerate}
We shall now determine the topological conditions for the existence of
an lpin structure.  According to Theorem \ref{t:strLpin}, one has to
distinguish two cases depending on the parity of the dimension of the
manifold.  If \(\dim M=m\) is even, then the group \(\Lpin_{m}\)
contains \( \Pin^{c}_{m} \) as a maximal compact subgroup and,
therefore, the existence of an lpin structure on an even-dimensional
manifold is equivalent to the existence of a pin\({}^c\) structure.
The subtler case of \( m \) odd is described in the following
\begin{thm}\label{t:exlpin}
A  manifold \( M\) of odd dimension admits an lpin structure if, and
only if, there exists  an element \( c\in {\mathrm
H}^2(M, {\mathbb Z})\)  and a
 real vector bundle \(E\) over \( M\), of fibre
dimension \( 2 \),  such that
\begin{equation}
w_2(TM)+w_2(E) \equiv c  \bmod 2.
\label{e:cond}
\end{equation}
\end{thm}

{\it Proof.} 
The set of all elements given in \neu{e:Aodd} generates the group
\( \Lpin_{m} \). Consider the homomorphism
\begin{equation}\label{e:homo}
\pi : \Lpin _m \rightarrow \Ort_{m} \times ({ \mathbb Z}_2
\times_{\varphi} \mathbb{C}^\times) \times \mathbb{C}^\times
\end{equation}
mapping an odd element \neu{e:Aodd} of \( \Lpin_{m} \) onto the
triple
\[
\bigl( \widetilde{\Ad}(a),
( - 1, \lambda\mu^{-1} ), \lambda\mu  \bigr).
\]
One easily checks that \( \pi \) provides a
two-fold covering. Consider the subgroups \( H_{0} \), \( H_1 \),  and
\( H_2 \) of   \( \Lpin _m \) defined by
\begin{eqnarray*}
H_{0} &=& \{  \mbox{\scriptsize\(\left(\begin{array}{cc}\sigma (a) &0\\  0
& \sigma (a) \end{array}\right)\)}\ts  \vert\ts  a \in \Pin _m \},
\\
H_1 &=& \{\mbox{\scriptsize\(\left(\begin{array}{cc}
\lambda &0\\0&\lambda^{-1}\end{array}\right),
\left(\begin{array}{cc} 0 & - {\lambda}^{-1} \\ \lambda &0
\end{array}\right)\)} \ts\vert \ts \lambda \in {\mathbb C}^\times
\}, \\
H_2 &=& \{ \mbox{\scriptsize \(\left(\begin{array}{cc}\lambda &0 \\  0 &
\lambda \end{array}\right)\)}\ts \vert\ts  \lambda \in
{\mathbb C}^\times \}.
\end{eqnarray*}
If \( i\neq j \), \( i,j=0,1,2 \), then \( H_{i}\cap H_{j}= {\mathbb
Z}_{2}\)
and every two elements of two different subgroups commute. Moreover,
the restrictions of the covering homomorphism \neu{e:homo}
to \( H_{0} \),  \( H_1 \) and \( H_2 \)  are two-fold coverings of
the groups
 \( \Ort_{m} \), \({\mathbb Z}_2 \times_{\varphi} {\mathbb C}^\times\)
and
\({\mathbb C}^\times\),  respectively.

\noindent An  lpin structure on the manifold \( M \) defines via the
covering \(\pi\):
\begin{itemize}
\item[\ (i)]  a two-dimensional real vector bundle \(E\) associated
with
the representation of the  group\\
\phantom{xxxxxx} \({\mathbb Z}_2 \times_{\varphi}
{\mathbb C}^{\times}= {\Ort}_{2} \times {\mathbb R}_+ \subset \GL_{2}(
{\mathbb R})\);
\item[(ii)] an oriented two-dimensional real vector bundle \(F\)
associated with the representation of the group\\
\phantom{xxxxxx} \( {\mathbb C}^{\times}
=\SO_{2}  \times {\mathbb R}_+ \subset \GL_{2}^+ ({\mathbb R})\).
\end{itemize}
Suppose that the bundles \(E\) and \(F\) are given. Fix a covering
\(\{
U_{\alpha} \}\) of the manifold \(M\) and denote by
\smallskip

\begin{tabular}{ll}
\(g_{\alpha \beta}:U_{\alpha} \cap U_{\beta} \to \Ort_m\) & the
transition functions of
 the \\ & tangent  bundle \(TM\),\\
\(h_{\alpha \beta} :U_{\alpha} \cap U_{\beta} \to {\mathbb Z}_2
\times_{\varphi} \mathbb{C}^{\times}\) & the transition functions of
the  \\
& bundle \(E\),\\
\(k_{\alpha \beta} :U_{\alpha} \cap U_{\beta} \to
\mathbb{C}^{\times}\)
&
the transition functions of the \\ & bundle \(F\).
\end{tabular}

\smallskip

\noindent Let \(g_{\alpha \beta}^*,  h_{\alpha \beta}^*, k_{\alpha
\beta}^*\) be their lifts into the covering groups:
\[
g_{\alpha \beta}^*, \, \, h_{\alpha \beta}^*, \, \, k_{\alpha
\beta}^* : U_{\alpha} \cap U_{\beta} \to  H_{0}, \, \, H_1, \, \,
H_2.
\]
We consider the map
\[
\varPhi_{\alpha \beta}^* = g_{\alpha \beta}^* h_{\alpha \beta}^*
k^*_{\alpha \beta} :U_{\alpha} \cap U_{\beta} \to \Lpin _m .
\]
The cocycle \(\varPhi^*_{\alpha \beta} \varPhi^*_{\beta \gamma}
\varPhi^{*-1}_{\alpha \gamma}\) represents in the \v{C}ech
cohomology 
group \({\mathrm H}^2 (M, {\mathbb Z}_2)\)
the obstruction to the reduction of the
structure group \(\Ort_m \times ({\mathbb Z}_2 \times_{\varphi}
\mathbb{C}^{\times} ) \times \mathbb{C}^{\times}\) to the structure
group  \(\Lpin
_m\) (see, e.g.,  [21]). Since pairs of  elements belonging to
different  subgroups
\(H_{0},  H_1, H_2\) commute, we obtain
\[
\varPhi^*_{\alpha \beta} \varPhi^*_{\beta \gamma}
\varPhi^{*-1}_{\alpha
\gamma} = \left( g^*_{\alpha \beta} g^*_{\beta \gamma}
g^{*-1}_{\alpha \gamma} \right) \left( h^*_{\alpha \beta} h^*_{\beta
\gamma} h^{*-1}_{\alpha \gamma} \right) \left(k^*_{\alpha \beta}
k^*_{\beta \gamma} k^{*-1}_{\alpha \gamma} \right) .
\]
The cocycle \( g^*_{\alpha \beta} g^*_{\beta \gamma} g^{*-1}_{\alpha
\gamma}\) represents the obstruction for the existence of a pin
structure on the manifold \(M\), i.e., the characteristic class
\(w_1(TM)^2 + w_2(TM)\). Similarly, the two other cocycles define
\(w_1 (E)^2 + w_2 (E)\) and \(w_2 (F)\), respectively.
Therefore, the vector bundle \(TM \oplus E \oplus F\) admits  a
reduction of the structure group to the group \(\Lpin _m\) if and
only if the condition
\begin{equation}\label{e:proof}
 w_1(TM)^2 + w_2(TM) + w_1(E)^2 + w_2 (E) + w_2
(F)=0
\end{equation}
in  holds.
Since \(F\) is an oriented vector
bundle its Stiefel-Whitney class \(w_2 (F)\) is the  mod 2
reduction of its Euler class \(e (F)\).
 Moreover, according to the Lemma, \(w_1(TM)^2\)  is also the  mod 2
reduction of some integral cohomology class. Therefore,
the existence of an lpin structure implies condition 
\neu{e:cond}. Conversely, suppose that $w_2 (TM) + w_2
(E)$ is the mod 2 reduction of some integral cohomology class. Using
the Lemma we conclude that $w_1 (TM)^2 + w_2 (TM) +
w_1 (E)^2 + w_2 (E)$ is also the mod 2 reduction of some
integral cohomology class $c \in H^2 (M, {\mathbb Z})$. There
exists an oriented vector bundle $F$
with fibres of real dimension two such that
its Euler class $e(F)$ coincides with the cohomology class $c \in H^2
(M, {\mathbb Z})$ (see, e.g., \cite{MS}). For this vector bundle $F$
equation \neu{e:proof} holds; therefore, $M$ admits an lpin structure. \qed \medskip 

\begin{exmp}
Let \(M\) be a manifold of dimension \( m=2n-1 \) isometrically
immersed in the Euclidean space \( {\mathbb R}^{2n+1} \).
 Since the codimension is two, there is  a natural
choice of a  two-dimensional bundle \(E\), namely the normal bundle
of \( M \).
 The Whitney theorem gives
\[
w_1(E)\ =\ w_1(TM),\quad w_1(TM)^2+w_2(TM)+w_2(E)=0.
\]
Therefore, every submanifold of codimension two of the Euclidean
space  admits an lpin structure. Moreover, in this case, there is an
explicit construction of the corresponding spinor bundle, similar to
the one known for hypersurfaces \cite{T}.
 Refer to Proposition \ref{p:CDP} and consider a
Pauli representation \( \theta:\Cl^{0}_{2n+1}\to \End S \) so that
\( S \) is complex \( 2^{n} \)-dimensional.  For every \(
x\in M \), the space \( T_{x}M \) can be identified with a \( (2n-1)
\)-dimensional vector subspace of \( {\mathbb R}^{2n+1} \).
Let \( \eta=e_{1}\dots e_{2n+1} \) be a volume element for
\( {\mathbb R}^{2n+1} \). Choose \( \iota\in\{1,\sqrt{-1}\} \) so that
\( \eta^{2}=\iota^{2} \). One makes
\( \varSigma =M\times S\to M \) into a spinor bundle by putting \(
\tau(v)(x,s)=(x,\iota\theta(v\eta)s )\) for \( v\in T_{x}M \subset
{\mathbb R}^{2n+1}\) and \( s\in S \).
\end{exmp}
\begin{exmp}
The latter example is a special case of a  more general situation
where a codimension two immersion induces an lpin structure on the
submanifold. Let \(N\)  be a  manifold of dimension \( 2n+1 \)
and assume that its
second Stiefel-Whitney class \(w_2 (TN)\) is the
 mod 2 reduction of some integral cohomology class \(c \in
{\mathrm H}^2 (N, {\mathbb Z})\).
For every codimension two submanifold \(M \) of \( N \) the formula
\[
w_2(TM) + w_2 (E) \equiv c \bmod 2
\]
holds, where \(E\) denotes again the two-dimensional normal bundle of
the submanifold \(M\). The condition on the manifold
\(N\) is satisfied, for example, in the following two cases:
\(N=X\times Y\), where \( X\) is a complex manifold  and  \( Y \) is
parallelizable;
\(N\) is a Sasakian manifold.
Therefore,  every submanifold  of codimension two in these spaces
admits an lpin structure.
\end{exmp}

\begin{exmp}\label{ex:grass}
Throughout this example \( M \) denotes the {\it Grassmann manifold\/}
\(G_{5,2}\)   of all (non-oriented)
2-dimensional linear subspaces of the 5-dimensional real vector space
\({\mathbb R}^5\).

The six-dimensional manifold \( M \) is  non-orientable,
connected and compact \cite{MS}. Its homology groups are  known
\cite{Th}:
\[
{\mathrm H}_1 (M,{\mathbb Z} )=
 {\mathbb Z}_2 , \qquad {\mathrm H}_2
(M,{\mathbb Z})= {\mathbb Z}_2.
\]
Using these results we find the \({\mathbb Z}\)- and the
\({\mathbb Z}_2\)-cohomology groups:
\begin{eqnarray*}
{\mathrm H}^1 (M,{\mathbb Z})=
0, \quad && \quad  {\mathrm H}^2 (M, {\mathbb Z})= {\mathbb Z}_2, \\
{\mathrm H}^1 (M,{\mathbb Z}_2 )=
{\mathbb Z}_2, \quad && \quad {\mathrm H}^2
(M,{\mathbb Z}_2 )= {\mathbb Z}_2 \oplus {\mathbb Z}_2.
\end{eqnarray*}
Denote by \(\gamma\) the canonical 2-dimensional vector bundle over
\(M\). Its first Stiefel-Whitney class \(w_1 (\gamma)\) is the
unique non-trivial element in \({\mathrm H}^1 (M,{\mathbb Z}_2)\):
\[
{\mathrm H}^1 (M,{\mathbb Z}_2 )=
 \{ 0, w_1 (\gamma) \}.
\]
Explicitly,
\[
{\mathrm H}^2 (M,{\mathbb Z}_2 )= \{ 0, w_1(\gamma)^2,
w_2 (\gamma), w_1 (\gamma)^2 + w_2 (\gamma) \}.
\]
Consider the restriction map \(r: {\mathbb Z}_2 =
{\mathrm H}^2 (M,
{\mathbb Z} ) \to {\mathrm H}^2 (M,{\mathbb Z}_2 )= {\mathbb Z}_2
\oplus {\mathbb Z}_2\). Since \(w_1(\gamma)^2 \) is  the
restriction of an integral cohomology class, the image of
\({\mathrm H}^2 (M,{\mathbb Z} )\) by \( r \) is
the set
\( \{ 0, w_1 (\gamma)^2 \} \).
The tangent bundle \(TM\) is isomorphic to the tensor product
 \(\gamma \otimes \gamma^{\perp}\), where \(\gamma^{\perp}\)
denotes the 3-dimensional vector bundle over \(M\) whose fibre at
the point (= plane in \( {\mathbb R}^{5} \)) \(\varPi \in M\)
is the 3-space \(\varPi ^{\perp}\). A
computation of the Stiefel-Whitney classes yields
\[
w_1 (TM )= w_1 (\gamma),\quad  w_2 (TM)= w_1 (\gamma)^2 + w_2
(\gamma).
\]

Consequently, the non-orientable manifold \(M\) does not admit a
 { pin}\(^c\)-structure. However, for the  bundle
\(E= \gamma\) over \(M\) we have
\[
w_2 (TM)+ w_2 (E) = w_1 (\gamma)^2
\]
and this class is the mod 2 reduction of an integral cohomology class
\(c \in {\mathrm H}^2 (M,{\mathbb Z})\).
 Consider now the
7-dimensional manifold
 \(M^\prime  =M \times {\mathbb R}\)
or \(M^\prime =M \times {\mathbb S}_{1}\).
Again, \(M^\prime\) is a non-orientable
manifold without a pin\(^c\)-structure. However, it satisfies the
condition of Theorem \ref {t:exlpin}; therefore, it admits an lpin
structure.
\end{exmp}

\vspace{1cm}

\begin{center}
{\scshape Acknowledgements}
\end{center}
The research reported in this paper was
supported in part by the Sonderforschungsbereich 288 at
Humboldt University, Berlin, 																	   %						   %
by the Polish Committee for Scientific Research (KBN) under grant
no. 2 P03B 017 12 at Warsaw University
and by the Foundation for Polish-German Cooperation
with funds provided by the Federal Republic of Germany.\\

% \end{article}
 \end{document}